\documentstyle[12pt, page1]{article}
\newcommand{\sect}[1]{\section{#1}\setcounter{equation}{0}}

\font\mbn=msbm10 scaled \magstep1
\font\mbs=msbm7 scaled \magstep1
\font\mbss=msbm5 scaled \magstep1
\newfam\mbff
\textfont\mbff=\mbn
\scriptfont\mbff=\mbs
\scriptscriptfont\mbff=\mbss
\def\mbf{\fam\mbff}

\def\Z{{\mbf Z}}
\def\Co{{\mbf C}}

\def\Di{{\mbf D}}
\newtheorem{Th}{Theorem}[section]
\newtheorem{Lm}[Th]{Lemma}
\newtheorem{C}[Th]{Corollary}
\newtheorem{D}[Th]{Definition}

\newtheorem{R}[Th]{Remark}
\newtheorem{Ex}[Th]{Example}
\author{Alexander Brudnyi\thanks
{Research supported in part by NSERC.
 }
\\ Department of Mathematics and Statistics\\ University of Calgary, 
Calgary\\ Canada} 
\title{Cyclicity and Maximal Multiplicity for Zeros of Families of Analytic
Functions}
\date{}
\begin{document}
\maketitle
\begin{abstract}
{Let $f_{\lambda}$ be a family of holomorphic functions in 
the unit disk $\Di\subset\Co$,
holomorphic in parameter $\lambda\in U\subset\Co^{n}$. 
We estimate the number of zeros of $f_{\lambda}$ in a
smaller disk via some characteristic of the ideal generated by
Taylor coefficients of $f_{\lambda}$. Our estimate is locally sharp and
improve the previous estimate obtained in [RY]. }
\end{abstract}
\sect{\hspace*{-1em}. Introduction.}
\noindent {\bf 1.1.}
In what follows $\Di_{r}:=\{z\in\Co\ :\ |z|<r\}$, \
$\overline{\Di}_{r}$ is the closure of $\Di_{r}$ and $\Di:=\Di_{1}$.
Let $U\subset\subset V\subset\Co^{n}$ be open connected sets, and
\begin{equation}\label{eq1}
f_{\lambda}(z)=\sum_{k=0}^{\infty}a_{k}(\lambda)z^{k},\ \ \ \ 
a_{k}\in {\cal O}(V),
\end{equation}
be a family of holomorphic functions in $\Di$
depending holomorphically on $\lambda\in V$.
Let ${\cal I}(f;U)$ be the 
ideal in ${\cal O}(U)$ generated by all $a_{k}(\lambda)$. 
Following the pioneering work of Bautin [B], we refer to
${\cal I}(f;U)$ as the {\em Bautin ideal} of $f_{\lambda}$ in $U$.
Further,   
$$
C(f;U):=\{\lambda\in U\ :\ f_{\lambda}\equiv 0\}
$$
is called the {\em central set} of $f_{\lambda}$ in $U$. The Hilbert 
finiteness theorem states that ${\cal I}(f;U)$ is generated by a finite 
number of coefficients. {\em The Bautin index} of $f$ in $U$ is 
the minimal number $d_{f}(U)$ such that $a_{0},...,a_{d_{f}(U)}$ generate 
${\cal I}(f;U)$. Usually computing ${\cal I}(f;U)$ and $d_{f}(U)$ 
is not easy. The number of zeros (counted with multiplicities) which 
$f_{\lambda}$ can have near $0$ for
$\lambda$ close to some $\lambda_{0}\in C(f;U)$ is called 
${\em cyclicity}$, following [R]. The next result was established by 
Yomdin [Y, Th. 3.1]:

{\em Assume that for any $\lambda\in V$ either 
$f_{\lambda}\equiv 0$, or the multiplicity of zero of $f_{\lambda}$ at
$0\in\Co$ is at most $N$. Let $I_{N}$ be the ideal in ${\cal O}(U)$
generated by the first $N$ Taylor coefficients $a_{0}(\lambda),\dots,
a_{N}(\lambda)$. Assume that $I_{N}$ is radical (i.e. $g^{s}\in I_{N}$
for some $s\geq 1$ implies that $g\in I_{N}$). Then the cyclicity of
$f_{\lambda}$, $\lambda\in U$, is at most $N$.}

This result follows from some theorems of [FY] based on the fact that 
$d_{f}(U)\leq N$. In particular, the required conditions are satisfied
if $f_{\lambda}$ depends linearly on $\lambda$. 
The main purpose of the present paper is to extend the above result
to a general situation. 

In [RY] it was shown that one can obtain a local 
upper bound on the number of zeros of 
$f_{\lambda}$ just in terms of $d_{f}(U)$. More precisely, 
there is some small positive $r<1$ depending on $d_{f}(U)$ and 
${\cal I}(f;U)$ such that each function $f_{\lambda}(z)$ has at most 
$d_{f}(U)$ complex zeros in the disk $\overline{\Di}_{r}$. In our paper
we will improve this estimate and will show that
another algebraic characteristic of ${\cal I}(f;U)$ is responsible
to the estimate of the number of zeros 
of $f_{\lambda}$. Moreover, our local estimate is sharp.\\
{\bf 1.2.} To formulate the results, let $K\subset\subset V$ be a compact. 
For any open $W$ by ${\cal O}_{c}(\Di; W)$ we denote the set of maps 
$\Di\mapsto W$ holomorphic in open neighbourhoods of $\overline{\Di}$. 
Assume that $W\subset V$. For any $\phi\in {\cal O}_{c}(\Di; W)$ consider 
the function
$$
f_{\phi(w)}(z)=\sum_{k=0}^{\infty}a_{k}(\phi(w))z^{k},\ \ \
w,z\in\Di\ .
$$
Let $d(\phi)$ be the Bautin index of the Bautin ideal of $f_{\phi(w)}$ in 
$\Di$.
Further, we fix a sequence of open sets $\{O_{j}\}$, 
$O_{j+1}\subset\subset O_{j}\subset\subset V$ for any $j$, such that 
$\cap_{j}O_{j}=K$.
\begin{D}\label{d1}
The integer number
$$
\mu_{f}(K):=
\lim_{j\to\infty}\sup_{\phi\in {\cal O}_{c}(\Di; O_{j})}d(\phi) 
$$
will be called the maximal multiplicity on $K$ of zero of $f_{\lambda}$ at
$0\in\Co$.
\end{D}
Obviously, $\mu_{f}(K)\leq d_{f}(O_{j})$, the Bautin 
index of ${\cal I}(f;O_{j})$ in $O_{j}$ (for any $j$). 
Therefore the definition is correct and $\mu_{f}(K)<\infty$.
Also, $\mu_{f}(K)$ does not depend of the choice of $\{O_{j}\}$. 

Below we give another characterization of 
$\mu_{f}(K)$. First we choose an open set $U$
such that $O_{j}\subset\subset U\subset\subset V$ for any $j$.
Let $N_{j}$ be minimum of integers $N$ for which there is 
$c(N)>0$ so that 
\begin{equation}\label{eq2}
|a_{k}(\lambda)|\leq c(N)\cdot\max_{U}|a_{k}|\cdot
\max_{i=0,...,N}|a_{i}(\lambda)|\ \ \ {\rm for}\ \ \ k>N,\
\lambda\in O_{j}\ .
\end{equation}
\begin{Th}\label{te1}
$$
\mu_{f}(K)=\lim_{j\to\infty}N_{j}\ .
$$
\end{Th}
We will show that one
can take the Bautin index $d_{f}(U)$ of ${\cal I}(f;U)$
as one of such $N$ in (\ref{eq2}). Let $c(N_{j})$ be the best constant in
(\ref{eq2}) for $N=N_{j}$. We set
\begin{equation}\label{best}
c_{\mu_{f}}(K)=\lim_{j\to\infty}c(N_{j})\ .
\end{equation}
The definition is correct because there is $i_{0}$ such that for any
$i\geq i_{0}$, $\displaystyle N_{i}=\lim_{j\to\infty}N_{j}$. Also,
$c_{\mu_{f}}(K)$ does not depend of the choice of $\{O_{i}\}$. In fact,
as follows from the proof of the theorem $c_{\mu_{f}}(K)$ {\em depends only 
on $K$ and $a_{0},\dots,a_{\mu_{f}}(K)$}. In general $c_{\mu_{f}}(K)$ cannot 
be estimated effectively. However, in many cases if the maximal multiplicity 
$\mu_{f}(K)$ is known, $c_{\mu_{f}}(K)$ can be found 
by a finite computation (which involves a resolution
of singularities type algorithm for the central set of $f_{\lambda}$).

Further, for $R<1$ and $\lambda\notin C(f;V)$, $\lambda\in V$, 
we set
\begin{equation}\label{eq3}
m_{f_{\lambda}}(R):=\sup_{\Di_{R}}\log |f_{\lambda}|\ \ \ \ 
{\rm and}\ \ \mu_{f}(\lambda,R):=m_{f_{\lambda}}(R)-m_{f_{\lambda}}(R/e)\ .
\end{equation}
\begin{Th}\label{te2}
Let
$$
\mu_{i}:=\limsup_{R\to 0}\sup_{\lambda\in O_{i}
\setminus C(f; O_{i})}\mu_{f}(\lambda,R).
$$
Then
$$
\mu_{f}(K)=\lim_{i\to\infty}\mu_{i}\ .
$$
\end{Th}
{\bf 1.3.} Next we formulate some simple corollaries from
Theorems \ref{te1} and \ref{te2}.
\\
(1)\ There is $i_{0}$ such that the central set $C(f;O_{i_{0}})$ is defined
as the set of common zeros of the first $\mu_{f}(K)+1$ Taylor coefficients
$a_{0}(\lambda),\dots,a_{\mu_{f}(K)}(\lambda)$ of $f_{\lambda}$
$(\lambda\in O_{i_{0}})$.\\
(2)\ Let $f_{\lambda}$ and $g_{\lambda}$ be families of holomorphic
functions in $\Di$, holomorphic in $\lambda\in V$,
and let $f', e^{f}$ be the families 
$f_{\lambda}'(z)=\frac{df_{\lambda}(z)}{dz}$ and $e^{f_{\lambda}}$, 
respectively. Then 
$$
\mu_{f}(K)\leq\mu_{f'}(K)+1\ ; \ \ \ \mu_{e^{f}}(K)=0\ ;\ \ \
\mu_{fg}(K)\leq\mu_{f}(K)+\mu_{g}(K)\ .
$$
(3)\ Let $S\subset V$ be a compact and ${\cal K}(S)$ be the space of all 
compact subsets of $S$ equipped with the Hausdorff metric. Then
the function $\mu_{f}:{\cal K}(S)\rightarrow\Z_{+}$ is
upper-semicontinuous, i.e., if $\{K_{i}\}$ is a sequence of compacts in $S$ 
converging in the Hausdorff metric to $K\subset S$, 
then 
$$
\limsup_{i\to\infty}\mu_{f}(K_{i})\leq\mu_{f}(K)\ .
$$
(4) For any family of compacts $\{K_{i}\}_{i\in I}$ in $V$ such that
$\cup_{i\in I}K_{i}$ is a compact 
$$
\mu_{f}(\cup_{i\in I}K_{i})=\max_{i\in I}\mu_{f}(K_{i})\ .
$$
(Here the maximum is taken because by (3) $\mu_{f}$ is upper-semicontinuous 
on ${\cal K}(\cup_{i}K_{i})$.) 
For $x\in V$ being a point the number $\mu_{f}(x)$ will 
be called {\em the generalized multiplicity} of zero of $f_{x}$ at 
$0\in\Co$. From 
Theorem \ref{te2} it follows that for $x\not\in C(f;V)$ the number 
$\mu_{f}(x)$ coincides with the usual multiplicity of zero of $f_{x}$ at 
$0\in\Co$. Also, for a compact $K$ from the above identity we have
$$
\mu_{f}(K)=\max_{x\in K}\mu_{f}(x)\ .
$$

We leave proofs of these simple properties as an exercise for the reader. \\
{\bf 1.4.} According to Theorem \ref{te1} there is $i_{0}$ such that 
$$
N_{i}=\mu_{f}(K)\ \ \ {\rm and}\ \ \ c(N_{i})\leq 2c_{\mu_{f}}(K)\ \ \
{\rm for\ any}\ \ \ i\geq i_{0}\ .
$$
Assume that 
\begin{equation}\label{eq5}
M:=\sup_{k>\mu_{f}(K),\lambda\in U}|a_{k}(\lambda)|<\infty\ 
\end{equation}
with $U$ as in Theorem \ref{te1}, and set
$$
R:=\frac{1}{4c_{\mu_{f}}(K)\cdot M\cdot 2^{30\mu_{f}(K)}+2}\ .
$$
Let $P_{\lambda}$ be the Taylor polynomial of $f_{\lambda}$
of degree $\mu_{f}(K)$. Let ${\cal N}_{r}(f_{\lambda})$ and
${\cal N}_{r}(P_{\lambda})$, $\lambda\in V$, be the number of 
zeros of $f_{\lambda}$ and $P_{\lambda}$ in $\overline{\Di}_{r}$.
\begin{Th}\label{te5}{\bf (Cyclicity Theorem)}
\begin{description}
\item[{\rm (1)}]\ \ \ ${\cal N}_{r/2}(P_{\lambda})\leq
{\cal N}_{r}(f_{\lambda})\leq {\cal N}_{2r}(P_{\lambda})$ for any \
$r<R$,\ $\lambda\in O_{i}$,\ $i\geq i_{0}$. 
\item[{\rm (2)}]\ \ For any $r< R$, $i\geq i_{0}$ there is some 
$\lambda\in O_{i}$ such that 
${\cal N}_{r}(f_{\lambda})=\mu_{f}(K)$.
\end{description}
\end{Th}
In particular, from (1) we have 
${\cal N}_{R}(f_{\lambda})\leq \mu_{f}(K)$ for any $\lambda\in O_{i}$, \
$i\geq i_{0}$.
\begin{R}\label{re1}
{\rm The straightforward application of Lemma 2.2.3 of [RY] 
gives also the following global estimate}
$$
{\cal N}_{1/4}(f_{\lambda})< 4\mu_{f}(K)+\log_{5/4}(2+2c_{\mu_{f}}(K)\cdot M)
\ \ \ {\rm for}\ \ \  \lambda\in O_{i},\ \ i\geq i_{0}\ .
$$
\end{R}
\begin{Ex}\label{e1}
{\rm (1) Let $B\subset\Co^{3}$ be a complex ball centered at 0 and}
$$
f_{\lambda}(z)=\lambda_{1}^{2}+\lambda_{2}^{2}z+\lambda_{3}^{2}z^{2}+
\lambda_{1}\lambda_{2}z^{4}+\lambda_{1}\lambda_{3}z^{5}+
\lambda_{2}\lambda_{3}z^{6},\ \ \ \ \ (\lambda,z)\in B\times\Di\ .
$$
{\rm It is easy to see that the Bautin ideal
in ${\cal O}(B)$ is generated by all coefficients of the function. Therefore
the Bautin index is 5 and according to [RY] the number of zeros of any 
$f_{\lambda}$ is $\leq 5$ in a small neighbourhood of $0\in\Co$. However, 
from the inequalities}
$$
|\lambda_{i}\lambda_{j}|\leq\max\{|\lambda_{1}|^{2},|\lambda_{2}|^{2},
|\lambda_{3}|^{2}\},\ \ \ \ 1\leq i,j\leq 3,
$$
{\rm it follows that $\mu_{f}(\overline{B})=2$. Thus according to the above 
theorems, the number of zeros of $f_{\lambda}$ in $\overline{\Di}_{r}$ with 
$r$ small enough is $\leq 2$.
Moreover, for any $\lambda=(0,0,\lambda_{3})\in B$, $\lambda_{3}\neq 0$, the
number of zeros (counted with multiplicities) of $f_{\lambda}$ in 
$\overline{\Di}_{r}$ is exactly 2.\\
(2) Let $f_{\lambda}(z)=\lambda (z^{10}-\lambda)$, $\lambda\in\Di$. In
this case the Bautin ideal is not radical and the result of Yomdin do not
apply (see [Y, page 363]). However, it is easy to see that 
$\mu_{f}(\lambda)=0$, $\lambda\neq 0$, and
$\mu_{f}(0)=\mu_{f}(\overline{\Di})=d_{f}(\Di)=10$. Thus the cyclicity is 10.
\\
(3) Assume that for any $\lambda\in V\subset\Co^{n}$ the function 
$f_{\lambda}$ satisfies 
\begin{equation}\label{dif}
f_{\lambda}^{(r)}(z)+p_{r-1\lambda}(z)f_{\lambda}^{(r-1)}(z)+\dots +
p_{1\lambda}(z)f_{\lambda}(z)=0,\ \ \ \ z\in\Di\ .
\end{equation}
Here each $p_{i\lambda}(z)$ is holomorphic in $(\lambda, z)\in V\times\Di$.
Then [RY, Corollary 4.2] and Theorem 
\ref{te1} above imply that $\mu_{f}(K)\leq r-1$ for any
compact $K\subset V$. Let 
$F_{\lambda}(z)=\sum_{k=1}^{m}P_{k}(z)e^{Q_{k}(z)}$
where $P_{k}, Q_{k}$ are holomorphic polynomials of maximal degree $p$ and 
$q$, respectively, and $\lambda\in\Co^{m(p+q+2)}$ is the vector of 
coefficients of all $P_{k},Q_{k}$ $(k=1,\dots,m)$. 
First consider $q=1$ (usual exponential polynomials).
Then $F_{\lambda}$ satisfies
(\ref{dif}) with $r=m(p+1)$. Hence $\mu_{F}(K)\leq m(p+1)-1$ for any
compact $K\subset\Co^{m(p+3)}$. Assume now that $q\geq 2$. 
It is easy to check (see e.g. [VPT]) that $F_{\lambda}$ satisfies
(\ref{dif}) with $r=\frac{(p+1)(q^{m}-1)}{q-1}$, and so
$\mu_{F}(K)\leq \frac{(p+1)(q^{m}-1)}{q-1}-1$. However, one obtains a 
better estimate using [Br, Lemma 8]. This result says that there
is $r_{0}>0$ such that for any $r\leq r_{0}$ the number of zeros in $\Di_{r}$
of any $F_{\lambda}$ is $\leq 3\cdot 2^{m-1}(p+q-1)$. Then
by Theorem \ref{te5},
$\mu_{F}(K)\leq 3\cdot 2^{m-1}(p+q-1)$.}
\end{Ex}
\sect{\hspace*{-1em}. Proofs of Theorems 1.2 and 1.3.}
\noindent 
{\bf 2.1. Resolution Theorem.} Our main tool is a version
of Hironaka's theorem on resolution of singularities proved in Theorem 4.4 
and Lemma 4.7 of Bierstone and Milman [BM]. As usual, 
if $z=(z_{1},\dots,z_{n})\in\Co^{n}$ and 
$\alpha=(\alpha_{1},\dots,\alpha_{n})$ is a multi-index 
($\alpha_{j}\in\Z_{+}$ ), $z^{\alpha}:=
z_{1}^{\alpha_{1}}\dots z_{n}^{\alpha_{n}}$ and 
$|\alpha|:=\sum_{1\leq j\leq n}\alpha_{j}$. If $\alpha$ and $\beta$ are
multi-indices we write $\alpha\prec\beta$ to mean that there is a multi-index
$\gamma$ such that $\beta=\alpha+\gamma$.

Fix $\lambda_{0}\in W$, an open subset
of $\Co^{n}$. A {\em dominating family} for $W$ at $\lambda_{0}$ is a
finite collection $(W_{\alpha},K_{\alpha},\phi_{\alpha})_{1\leq\alpha\leq A}$
where, for each $\alpha$, $K_{\alpha}\subset\subset W_{\alpha}$, 
$W_{\alpha}$ is an open set in $\Co^{n}$ containing $0$, 
$\phi_{\alpha}: W_{\alpha}\longrightarrow W$ is a holomorphic map satisfying 
the two conditions
\begin{description}
\item[{\rm (1)}]\ \ 
$det\phi_{\alpha}'\neq 0$ outside a complex analytic variety
of codimension $\geq 1$, and
\item[{\rm (2)}]
\ \ The images $\phi_{\alpha}(K_{\alpha})$ $(\alpha=1,\dots,A)$
cover a neighbourhood of $\lambda_{0}$ in $\Co^{n}$.
\end{description}
{\bf BM Theorem.} {\em  Let $f_{1},\dots,f_{N}$ be holomorphic functions
defined on a neighbourhood of $\lambda_{0}\in\Co^{n}$. Suppose that
none of the $f_{j}$ vanishes identically in any neighbourhood of 
$\lambda_{0}$. Then there exists a dominating family 
$(W_{\alpha},K_{\alpha},\phi_{\alpha})_{1\leq\alpha\leq A}$ for $W$ at
$\lambda_{0}$, such that for each $\alpha$ we can find multi-indices
$\gamma_{1\alpha},\dots,\gamma_{N\alpha}$ and functions
$h_{1\alpha},\dots,h_{N\alpha}$ on $W_{\alpha}$ with the following
properties
\begin{description}
\item[{\rm (A)}]\  Each $h_{j\alpha}$ is a nowhere-vanishing holomorphic 
function on $W_{\alpha}$.
\item[{\rm (B)}]\ 
$f_{j}\circ\phi_{\alpha}(z)=h_{j\alpha}(z)\cdot z^{\gamma_{j\alpha}}$
for all $z\in W_{\alpha}$, $1\leq j\leq N$, $1\leq\alpha\leq A$.
\item[{\rm (C)}]\ 
For each $\alpha$, the multi-indices $\gamma_{1\alpha},\dots,
\gamma_{N\alpha}$ are totally ordered under $\prec$ (i.e. given
$1\leq i,j\leq N$, we either have $\gamma_{i\alpha}\prec\gamma_{j\alpha}$
or $\gamma_{j\alpha}\prec\gamma_{i\alpha}$).
\end{description} }
{\bf 2.2. Proof of Theorem \ref{te1}.} Let
$f_{\lambda}(z)=\sum_{k=0}^{\infty}a_{k}(\lambda)z^{k},\
a_{k}\in {\cal O}(V)$, and $K\subset\subset U\subset\subset V$ where
$K$ is compact and $U$ is open. Let $d_{f}(U)$ be the Bautin index of
${\cal I}(f;U)$ in $U$. Let $\lambda_{0}\in U$ and $W\subset U$ be an open
neighbourhood of $\lambda_{0}$. Assume that 
$(W_{\alpha},K_{\alpha},\phi_{\alpha})_{1\leq\alpha\leq A}$ is the
dominating family for $W$ at $\lambda_{0}$ for which 
$a_{0},\dots,a_{d_{f}(U)}$ satisfy BM Theorem. Here
$a_{j}\circ\phi_{\alpha}(z)=h_{j\alpha}(z)\cdot z^{\gamma_{j\alpha}}$
for all $z\in W_{\alpha}$, $1\leq j\leq d_{f}(U)$, 
$1\leq\alpha\leq A$. Then from (C) it follows that one can find an index 
$j_{0}(\alpha)$ and multi-indices $\tilde\gamma_{j\alpha}$ so that
$\gamma_{j\alpha}=\gamma_{j_{0}(\alpha)\alpha}+\tilde\gamma_{j\alpha}$, 
with $\gamma_{j_{0}(\alpha)\alpha}=0$. The minimal 
$j_{0}(\alpha)$ satisfying this condition will be denoted  
$j(W_{\alpha})$.
\begin{Lm}\label{le1}
For any $k\geq 0$ the function $\displaystyle 
g_{k\alpha}(z):=\frac{a_{k}\circ\phi_{\alpha}(z)}{
z^{\gamma_{j(W_{\alpha})\alpha}}}$
is holomorphic on $W_{\alpha}$. There is a constant 
$C(W_{\alpha})>0$ such that 
\begin{equation}\label{eq6}
|a_{k}\circ\phi_{\alpha}(z)|\leq C(W_{\alpha})\cdot\max_{W}
|a_{k}|\cdot\max_{i=0,\dots,j(W_{\alpha})}|a_{i}\circ\phi_{\alpha}(z)|\ ,
\ \ \ z\in K_{\alpha} .
\end{equation}
\end{Lm}
{\bf Proof.} For the first statement it suffices to consider 
$k>d_{f}(U)$. Then $a_{k}$ belongs to the ideal generated by
$a_{0},\dots,a_{d_{f}(U)}$ on $U$. So there are 
$b_{j}\in {\cal O}(U)$, $1\leq j\leq d_{f}(U)$, such 
that $a_{k}=\sum_{j=1}^{d_{f}(U)}b_{j}\cdot a_{j}$.
This implies the first part of the lemma. Further, let us cover the compact
$\overline{K}_{\alpha}$ 
by a finite number of closed polydisks $\Delta_{l}:=\{z\in\Co^{n}\ : 
\max_{1\leq i\leq n}|z_{i}-\xi_{il}|\leq a_{l}\}$, 
$\Delta_{l}\subset W_{\alpha}$, 
$\min_{1\leq i\leq n}||\xi_{il}|-|a_{l}||:=r_{l}>0$ $(l=1,\dots,s)$.
Assume that 
$z^{\gamma_{j(W_{\alpha})\alpha}}=z_{1}^{\gamma_{1}}\dots z_{n}^{\gamma_{n}}$ 
with $|\gamma_{j(W_{\alpha})\alpha}|=\sum_{i=1}^{n}\gamma_{i}$.
Then on the boundary torus $S_{l}^{n}\subset \Delta_{l}$ we have
$$
\min_{S_{l}^{n}}|z^{\gamma_{j(W_{\alpha})\alpha}}|\geq 
r_{l}^{|\gamma_{j(W_{\alpha})\alpha}|}:=
r_{l}(W_{\alpha})\ .
$$
Also, for $z\in\Delta_{l}$  by the definition we have
$$
|z^{\gamma_{j(W_{\alpha})\alpha}}|\leq
|a_{j(W_{\alpha})}\circ\phi_{\alpha}(z)|\cdot\max_{S_{l}^{n}}\frac{1}{
|h_{j(W_{\alpha})\alpha}|}\leq k_{l}(W_{\alpha})\cdot
\max_{i=0,\dots,j(W_{\alpha})}
|a_{i}\circ\phi_{\alpha}(z)|\ .
$$
Combining these inequalities we obtain (for $z\in\Delta_{l}$)
$$
\begin{array}{c}
\displaystyle
|a_{k}\circ\phi_{\alpha}(z)|=|g_{k\alpha}(z)|\cdot 
|z^{\gamma_{j(W_{\alpha})\alpha}}|
\leq k_{l}(W_{\alpha})\cdot
\max_{S_{l}^{n}}|g_{k\alpha}|\cdot\max_{i=0,\dots,j(W_{\alpha})}
|a_{i}\circ\phi_{\alpha}(z)|\leq\\
\displaystyle
k_{l}(W_{\alpha})\cdot\max_{S_{l}^{n}}|a_{k}\circ\phi_{\alpha}(z)|\cdot
\max_{S_{l}^{n}}\frac{1}{|z^{\gamma_{j(W_{\alpha})\alpha}}|}\cdot
\max_{i=0,\dots,j(W_{\alpha})}|a_{i}\circ\phi_{\alpha}(z)|\leq\\
\displaystyle
\frac{k_{l}(W_{\alpha})}{r_{l}(W_{\alpha})}\cdot\max_{W}|a_{k}|
\cdot\max_{i=0,\dots,j(W_{\alpha})}|a_{i}\circ\phi_{\alpha}(z)|\ .
\end{array}
$$
From here it follows the required inequality (\ref{eq6}) with
$$
C(W_{\alpha}):=\max_{1\leq l\leq s}
\frac{k_{l}(W_{\alpha})}{r_{l}(W_{\alpha})}\ .
\ \ \ \ \ \Box
$$
Set $j(W)=\max_{\alpha}j(W_{\alpha})$. Consider the family
$\{O_{i}\}$ as in Theorem \ref{te1}, i.e., $K\subset O_{i}$,
$O_{i+1}\subset\subset O_{i}\subset\subset U$ for any $i$. Let
$i_{0}$ be such that $\displaystyle \lim_{j\to\infty}N_{j}=N_{i}$ 
for any $i\geq i_{0}$.
\begin{Lm}\label{le2}
For any $i\geq i_{0}$ there are finite number
of points $\lambda_{li}\in\overline{O}_{i+1}$, open sets 
$W_{li}\subset O_{i}$, $\lambda_{li}\in W_{li}$,
and dominating families
$(W_{\alpha, li},K_{\alpha, li},\phi_{\alpha, li})_{1\leq\alpha\leq A_{li}}$
for $W_{li}$ at $\lambda_{li}$, satisfying BM Theorem for the functions
$a_{0},\dots,a_{d_{f}(U)}$ ($1\leq l\leq t_{i})$,
such that 
$$
\overline{O}_{i+1}\subset
\bigcup_{l=1}^{t_{i}}\bigcup_{\alpha=1}^{A_{li}}
\phi_{\alpha,li}(K_{\alpha,li})\ .
$$
\end{Lm}
{\bf Proof.} The proof follows from the fact that the closure 
$\overline{O}_{i+1}$ is compact. \ \ \ \ \ $\Box$

Let ${\cal W}_{i}=(W_{li})_{l=1}^{t_{i}}$ be the corresponding open cover of
$\overline{O}_{i+1}$. Then we set 
$$
U({\cal W}_{i})=\bigcup_{l=1}^{t_{i}}W_{li}\ \ \ {\rm and}\ \ \
j({\cal W}_{i})=\max_{1\leq l\leq t_{i}}j(W_{li})\ .
$$
Clearly, $U({\cal W}_{i})\subset O_{i}$, 
$U({\cal W}_{i+1})\subset\subset U({\cal W}_{i})$, and
$j({\cal W}_{i})\leq d_{f}(U)$. Now inequality (\ref{eq6}) applied 
to the elements of dominating families of a cover ${\cal W}_{i}$ as above
implies
\begin{Lm}\label{le3}
There is $C_{i}>0$ such that
$$
|a_{k}(\lambda)|\leq C_{i}\cdot\max_{U}|a_{j}|\cdot
\max_{i=0,...,j({\cal W}_{i})}|a_{i}(\lambda)| \ \ \ {\rm for}\ \ \
k\geq 0,\ \lambda\in O_{i+1}\ .\ \ \ \ \ \Box
$$
\end{Lm}
Next we will prove 
\begin{Lm}\label{le4}
For any $i\geq i_{0}$,
$$
\lim_{j\to\infty}N_{j}=j({\cal W}_{i}).
$$
\end{Lm}
{\bf Proof.} By definition,
$O_{i+1}\subset\subset O_{i_{0}}$.
Thus from Lemma \ref{le3} and from the definition of $N_{i+1}$ we have 
$\displaystyle L:=\lim_{j\to\infty}N_{j}\leq j({\cal W}_{i})$. Assume, to the 
contrary, that $L<j({\cal W}_{i})$. By definition, there is an element 
$(W_{\alpha},K_{\alpha},\phi_{\alpha})$ of one of the dominating families for
${\cal W}_{i}$ (as in Lemma \ref{le2}) such that 
$j(W_{\alpha})=j({\cal W}_{i})$. Since $U({\cal W}_{i})\subset\subset O_{i}
\subset\subset O_{i_{0}}$, we have from our assumption
\begin{equation}\label{eq7}
|a_{j({\cal W}_{i})}\circ\phi_{\alpha}(z)|\leq c(L)
\cdot\max_{U}|a_{j({\cal W}_{i})}|\cdot
\max_{j=0,...,L}|a_{j}\circ\phi_{\alpha}(z)|\ ,\ \ \ z\in W_{\alpha}\ .
\end{equation}
But by BM Theorem,
$a_{j}\circ\phi_{\alpha}(z)=h_{j\alpha}(z)\cdot z^{\gamma_{j\alpha}}$ 
where each $h_{j\alpha}$ 
is nowhere-vanishing on $W_{\alpha}$, and $j({\cal W}_{i})$ is the 
minimal number such that $\gamma_{j({\cal W}_{i})\alpha}\prec\gamma_{j\alpha}$
$(j=1,\dots,d_{f}(U))$. Then
(\ref{eq7}) gives a contradiction with minimality of $j({\cal W}_{i})$ 
(since $0\in W_{\alpha}$).
\ \ \ \ \ $\Box$

We set $j(K):=j({\cal W}_{i})$ $(i\geq i_{0})$.
It remains to prove that $\mu_{f}(K)=j(K)$. Assume also that the above 
$i_{0}$ is so big that for any $i\geq i_{0}$,
$$
\mu_{f}(K)=\sup_{\phi\in {\cal O}_{c}(\Di; O_{i})}d(\phi)\ .
$$
Let $\phi\in {\cal O}_{c}(\Di;O_{i_{0}})$ be such that $d(\phi)=\mu_{f}(K)$.
By definition $\phi$ maps some $\Di_{r}$, $r>1$, into $O_{i_{0}}$. Then 
the definition of $N_{i_{0}}=j(K)$ implies that
$$
|a_{k}\circ\phi (z)|\leq c(j(K))
\cdot\max_{U}|a_{k}|\cdot
\max_{j=0,...,j(K)}|a_{j}\circ\phi(z)|\ ,\ \ \ k>j(K),\ z\in\Di_{r}\ .
$$ 
Let $Z\subset\Di_{r}$ be the  set of common zeros of $a_{j}\circ\phi$, 
$0\leq j\leq j(K)$, counted with multiplicities, and 
$B_{Z}$ be the Blashke product in $\Di_{r}$ whose set of zeros is $Z$.
For any $i\geq 0$ we set $h_{i}:=\frac{a_{i}\circ\phi}{B_{Z}}$. From the
above inequality it follows that $h_{i}$ is holomorphic on $\Di_{r}$. Then
for some $r_{1}$, $1<r_{1}<r$,
$$
\max_{j=0,...,j(K)}|h_{j}(z)|\geq C>0,\ \ \
z\in\Di_{r_{1}}\ .
$$
Hence, by the corona theorem, there are bounded holomorphic on
$\Di_{r_{1}}$ functions $g_{0},\dots, g_{j(K)}$ such that
$$
\sum_{i=0}^{j(K)}g_{i}\cdot h_{j}\equiv 1\ \ \ {\rm on}\ \ \
\Di_{r_{1}}\ .
$$
From here it follows (for any $k$)
$$
a_{k}\circ\phi\equiv\sum_{i=0}^{j(K)}(a_{i}\circ\phi)\cdot 
(h_{k}\cdot g_{i})
\ \ \ {\rm on}\ \ \ \Di\ .
$$
This means that the Bautin index $d(\phi)\ (=\mu_{f}(K))$ is $\leq$ $j(K)$.

Conversely, let ${\cal W}_{i}$ and $W_{\alpha},K_{\alpha},\phi_{\alpha}$
be the same as in the proof of Lemma \ref{le4}. Here
$j(W_{\alpha})=j(K)$. Then from BM Theorem we have
$a_{j}\circ\phi_{\alpha}(z)=h_{j\alpha}(z)\cdot z^{\gamma_{j\alpha}}$, 
$z\in W_{\alpha}$, with a nowhere-vanishing $h_{j\alpha}$
$(1\leq j\leq d_{f}(U))$. Let us consider the closed disk 
$D_{s}:=\{(z,\dots,z)\in\Co^{n}\ :\ z\in\Co,\ |z|\leq s\}$ with $s$ so
small that $D_{s}\subset W_{\alpha}$. Then by the definition of $j(K)$, the
multiplicity of zero of each $(a_{k}\circ\phi_{\alpha})|_{D_{s}}$, 
$0\leq k<j(K)$, at $(0,\dots,0)\in D_{s}$ is greater than 
$|\gamma_{j(K)\alpha}|$, but
the multiplicity of zero of $(a_{j(K)}\circ\phi_{\alpha})|_{D_{s}}$ at
$(0,\dots,0)\in D_{s}$ equals $|\gamma_{j(K)\alpha}|$.
Now, for some $r>1$ we still have
$\tau(w):=(sw,\dots,sw)\in W_{\alpha}$ for $w\in\Di_{r}$. Let us define
$\phi\in {\cal O}_{c}(\Di;O_{i_{0}})$ as $\phi:=\phi_{\alpha}\circ\tau$.
According to the above argument for multiplicities, the Bautin index 
$d(\phi)$ of $f_{\phi(w)}$ in $\Di$ is $\geq j(K)$. This shows that 
$\mu_{f}(K)=j(K)$ .

The proof of the theorem is complete.\ \ \ \ \ \ $\Box$\\
{\bf 2.3. Proof of Theorem \ref{te2}.} Let $i_{0}$ be the same as in
the proof of Theorem \ref{te1}. We will prove that {\em for any 
$i\geq i_{0}+1$, $\mu_{i}=\mu_{f}(K)$.}

Let ${\cal W}_{i-1}$ be a cover from Lemma \ref{le2}.
Here $U({\cal W}_{i-1})\subset O_{i-1}\subset O_{i_{0}}$ and 
$j({\cal W}_{i-1})=\mu_{f}(K)$. By definition
there are a sequence of points $w_{s}\in O_{i}\setminus 
C(f;O_{i})$ and numbers $R_{s}>0$, $\displaystyle
\lim_{s\to\infty}R_{s}=0$, such that
$\displaystyle \mu_{i}=\lim_{s\to\infty}\mu_{f}(w_{s},R_{s})$.
Without loss of generality we may assume that $\displaystyle 
\lim_{s\to\infty}w_{s}=w\in \overline{O}_{i}$. Then from Lemma \ref{le2}
it follows that there is a dominating family
$(W_{\alpha},K_{\alpha},\phi_{\alpha})_{1\leq\alpha\leq A}$ for one of
the open sets of the cover ${\cal W}_{i-1}$ such that  images
$\phi_{\alpha}(\overline{K}_{\alpha})$ $(1\leq\alpha\leq A)$ cover a
neighbourhood of $w$. In particular,  we can find some $\alpha$ and a 
sequence $\{\widetilde w_{k}\}\subset K_{\alpha}$ such that 
$\phi_{\alpha}(\widetilde w_{k})=w_{s_{k}}$ for some subsequence 
$\{w_{s_{k}}\}\subset\{w_{s}\}$ and 
$\displaystyle \lim_{k\to\infty}\widetilde w_{k}:=\widetilde w\in W_{\alpha}$.
Now, according to BM Theorem and Lemma \ref{le1} the function
$$
h_{z}(y):=
\frac{f_{\phi_{\alpha}(z)}(y)}{z^{\gamma_{j(W_{\alpha})\alpha}}},\ 
\ \ z\in W_{\alpha},\ y\in\Di,
$$
is holomorphic and its central set in $W_{\alpha}$ is empty. 
Let $\mu$ be the multiplicity of zero of $h_{\widetilde w}$ at 0. Then the
function $g_{z}(y):=\frac{h_{z}(y)}{y^{\mu}}$ is nowhere-vanishing in a small
neighbourhood $O$ of $(\widetilde w,0)\in W_{\alpha}\times\Di$. Without
loss of generality we may assume that all pairs 
$(\widetilde w_{k},y)$, $y\in\overline{\Di}_{R_{s_{k}}}$, belong to $O$. 
Then we have
\begin{equation}\label{mu}
\mu_{i}=\lim_{k\to\infty}\mu_{f}(w_{s_{k}},R_{s_{k}})=
\lim_{k\to\infty}\mu_{h}(\widetilde w_{k},R_{s_{k}})=
\lim_{k\to\infty}\mu_{g}(\widetilde w_{k},R_{s_{k}})+\mu=\mu\ .
\end{equation}
But in $W_{\alpha}$ the maximal multiplicity of zero of $h_{z}$ at 0 
is $j(W_{\alpha})$ because by BM Theorem the Taylor coefficient of $h_{z}$ 
whose number is $j(W_{\alpha})+1$ is nowhere-vanishing and the previous
coefficients have common zero at $0\in W_{\alpha}$. This shows that
$$
\mu_{i}\leq j({\cal W}_{i-1})=\mu_{f}(K)\ .
$$

Let us prove the opposite inequality.
Recall that from Theorem \ref{te1} applied to the cover ${\cal W}_{i}$ 
it follows that there is an element 
$(W_{\alpha},K_{\alpha},\phi_{\alpha})$ of a dominating family of
the cover ${\cal W}_{i}$ (with $U({\cal W}_{i})\subset O_{i}$)
such that $j(W_{\alpha})=\mu_{f}(K)$. In particular, there is a point
$z_{0}\in W_{\alpha}$ such that the multiplicity $\mu$ of zero of 
$h_{z_{0}}$ (defined as above) at $0$ is $j(W_{\alpha})$. Then from 
(\ref{mu}) with $\displaystyle \lim_{k\to\infty}\widetilde w_{k}=z_{0}$
and from maximality of $\mu_{i}$ we have 
$$
\mu_{i}\geq\mu=j(W_{\alpha})=\mu_{f}(K)\ .
$$
\\
Combining these inequalities we get
$$
\mu_{i}=\mu_{f}(K)\ \ \ {\rm and}\ \ \ \lim_{i\to\infty}\mu_{i}=\mu_{f}(K)
$$
which completes the proof of the theorem.\ \ \ \ \ $\Box$
\sect{\hspace*{-1em}. Proof of the Cyclicity Theorem.}
\noindent
{\bf 3.1. Cartan's Lemma.} In the proof we use a version of the Cartan Lemma
proved in Levin's book [L, p.21].\\
{\bf Cartan's Lemma.} {\em Let $f(z)$ be a holomorphic function on 
$\Di_{2eR}$,
$f(0)=1$, and $\eta$ be a positive number $\leq\frac{3e}{2}$. Then there is 
a set of disks $\{D_{j}\}$ with $\sum_{j}r_{j}\leq 4\eta R$, where $r_{j}$
is radius of $D_{j}$ such that
$$
\log|f(z)|>-H(\eta)\log\max_{\Di_{2eR}}|f|
$$
for any $z\in\Di_{R}\setminus(\cup_{j}D_{j})$. Here $H(\eta)=2+
\log\frac{3e}{2\eta}$.}

Let $g(z)$ be a holomorphic function on $\overline{\Di}_{(6e+1)r/2}$. Let
$m_{1}:=\max_{\Di_{r/2}}|g|$ and $m_{2}:=\max_{\Di_{(6e+1)r/2}}|g|$. In what
follows $S_{t}:=\{z\in\Co\ :\ |z|=t\}$. From Cartan's Lemma we have
\begin{Lm}\label{le31}
There is a number $t_{r}$, \ $r/2\leq t_{r}\leq r$, such that
$$
\min_{S_{t_{r}}}\log |g|>\log m_{1}-7\cdot\log\frac{m_{2}}{m_{1}}\ .
$$
\end{Lm}
{\bf Proof.} Let $w\in S_{r/2}$ be such that $|g(w)|=m_{1}$.
We set 
$f(z):=\frac{g(z+w)}{g(w)}$.
Then $f$ is defined on $\Di_{2eR}$ with $R=\frac{3r}{2}$, and $f(0)=1$,
$\max_{\Di_{2eR}}|f|\leq\frac{m_{2}}{m_{1}}$. For any $t$, \
$r/2\leq t\leq r$, we set
$$
K_{t}:=\{z\in S_{t}\ :\ |f(z)|=\min_{S_{t}}|f|\}\ \ \ \ 
{\rm and}\ \ \ \ K=\cup_{t}K_{t}\ .
$$
Clearly one cannot cover $K$ by a set of disks $\{D_{j}\}$ with
$\sum_{j}r_{j}<4\eta R$ where $\eta=1/24$. In particular, by Cartan's lemma
there is $t_{r}\in [r/2,r]$ such that
$$
\min_{S_{t_{r}}}\log |f|\geq -(2+\log 36e)\log\frac{m_{2}}{m_{1}}>
-7\log\frac{m_{2}}{m_{1}}\ .
$$
Going back to $g$ gives the required inequality.\ \ \ \ \ $\Box$

Assume that $g$ as above is a polynomial of degree $d$. Then we have
\begin{C}\label{cor1}
There is a number $t_{r}$,\ $r/2\leq t_{r}\leq r$, such that
$$
\min_{S_{t_{r}}}|g|>\frac{1}{(6e+1)^{7d}}\cdot\max_{\Di_{r/2}}|g|>
\frac{1}{2^{29d}}\cdot\max_{\Di_{r/2}}|g|\ .
$$
\end{C}
{\bf Proof.} The Bernstein Doubling inequality for polynomials implies
$\frac{m_{2}}{m_{1}}\leq (6e+1)^{d}$.
Then we apply Lemma \ref{le31}.\ \ \ \ \ $\Box$\\
{\bf 3.2. Proof of the Cyclicity Theorem.} (1) The definition of $i_{0}$,
estimate (\ref{eq5}) and Theorem \ref{te1} imply for $i\geq i_{0}$,
$$
|a_{k}(\lambda)|\leq 2c_{\mu_{f}}(K)\cdot M\cdot \max_{j=0,...,\mu_{f}(K)}
|a_{j}(\lambda)|,\ \ \ \ \ k>\mu_{f}(K),\ \ \lambda\in O_{i}\ .
$$
Then in a disk $\Di_{R}$ with $R<1$ for $\lambda\in O_{i}$ we have
\begin{equation}\label{eq31}
|f_{\lambda}(z)-P_{\lambda}(z)|=
|\!\sum_{i>\mu_{f}(K)}a_{i}(\lambda)\cdot z^{i}|\leq 
2c_{\mu_{f}}(K)\cdot M\cdot
\frac{R^{\mu_{f}(K)+1}}{1-R}\cdot\max_{j=0,...,\mu_{f}(K)}|a_{j}(\lambda)|\ .
\end{equation}
Also, by the Cauchy inequality we have
$$
\max_{j=0,...,\mu_{f}(K)}|a_{j}(\lambda)|\leq\max_{\Di}|P_{\lambda}|\ .
$$
The last inequality, Corollary \ref{cor1} and the Bernstein Doubling 
inequality imply that there is $t_{R}$, \ $R/2\leq t_{R}\leq R$, such that
\begin{equation}\label{eq32}
\min_{S_{t_{R}}}|P_{\lambda}|>
\frac{\max_{\Di_{R/2}}|P_{\lambda}|}{2^{29\mu_{f}(K)}}\geq
\frac{R^{\mu_{f}(K)}\cdot\max_{\Di}|P_{\lambda}|}{2^{30\mu_{f}(K)}}\geq
\frac{R^{\mu_{f}(K)}}{2^{30\mu_{f}(K)}}\cdot
\max_{j=0,...,\mu_{f}(K)}|a_{j}(\lambda)|\ .
\end{equation}
We set
$$
R_{0}=
\frac{1}{2c_{\mu_{f}}(K)\cdot M\cdot 2^{30\mu_{f}(K)}+1}\ .
$$
Combining (\ref{eq31}) and (\ref{eq32}) we have for any $R<R_{0}$, 
$z\in S_{t_{R}}$, $\lambda\in O_{i}$,
$$
\begin{array}{c}
\displaystyle
|f_{\lambda}(z)-P_{\lambda}(z)|\leq 2c_{\mu_{f}}(K)\cdot M\cdot
\frac{R^{\mu_{f}(K)+1}}{1-R}\cdot\max_{j=0,...,\mu_{f}(K)}|a_{j}(\lambda)|<\\
\\
\displaystyle
\frac{R^{\mu_{f}(K)}}{2^{30\mu_{f}(K)}}\cdot
\max_{j=0,...,\mu_{f}(K)}|a_{j}(\lambda)|<\min_{S_{t_{R}}}|P_{\lambda}|\leq
|P_{\lambda}(z)|\ .
\end{array}
$$
From here by the Rouch\'{e} theorem it follows that 
$f_{\lambda}$ and $P_{\lambda}$ have the same number of zeros in 
$\overline{\Di}_{t_{R}}$.
If we apply the last statement to any $R<R_{0}/2$ we obtain
$$
{\cal N}_{R/2}(P_{\lambda})\leq {\cal N}_{t_{R}}(P_{\lambda})=
{\cal N}_{t_{R}}(f_{\lambda})\leq
{\cal N}_{R}(f_{\lambda})\leq {\cal N}_{t_{2R}}(f_{\lambda})=
{\cal N}_{t_{2R}}(P_{\lambda})\leq {\cal N}_{2R}(P_{\lambda})\ .
$$
This proves the first part of the theorem.\\
(2) Let $r<R_{0}/2$ and $i\geq i_{0}$ with $R_{0}$ and $i_{0}$ as above.
From the proof of Theorem \ref{te2} we know that
there is a point $\lambda\in O_{i}$, an open set $W\subset O_{i}$,
$\lambda\in W$, and an element $(W_{\alpha},K_{\alpha},\phi_{\alpha})$
of the dominating family for $W$ at $\lambda$ such that 
$j(W_{\alpha})=\mu_{f}(K)$. Moreover from (\ref{mu}) it follows that
for the function
$$
h_{z}(y):=
\frac{f_{\phi_{\alpha}(z)}(y)}{z^{\gamma_{j(W_{\alpha})\alpha}}},\ 
\ \ z\in W_{\alpha},\ y\in\Di\ , 
$$
there is $z_{0}\in W_{\alpha}$ such that the multiplicity of zero of
$h_{z_{0}}$ at $0$ equals $\mu_{f}(K)$. Further, we can find $r'$,
$r\leq r'<R_{0}/2$, such that $h_{z_{0}}|_{S_{r'}}$ is nowhere-vanishing.
In particular, there is an open connected neighbourhood 
$O\subset W_{\alpha}$ of 
$z_{0}$ such that for any $z\in O$ we still have that $h_{z}$ is 
nowhere-vanishing on $S_{r'}$. For $z\in O$ consider the integral
$$
{\cal N}_{r'}(h_{z})=\frac{1}{2\pi i}\int_{S_{r'}}
\frac{\frac{\partial}{\partial y}h_{z}(y)}{h_{z}(y)}dy
$$
which counts the number of zeros of $h_{z}$ in $\overline{\Di}_{r'}$.
Clearly the function ${\cal N}_{r'}(h_{z})$, $z\in O$, is continuous and 
integer-valued. Thus it is constant on $O$. From here by the
definition of the dominating family it follows that
there is a point $w\in O$ such that
$\lambda_{0}:=\phi_{\alpha}(w)\in O_{i}\setminus C(f;O_{i})$ and
$$
{\cal N}_{r'}(f_{\lambda_{0}})={\cal N}_{r'}(h_{w})={\cal N}_{r'}(h_{z_{0}})
\geq\mu_{f}(K)\ .
$$
But according to part (1) of the Cyclicity Theorem, for any 
$\lambda\in O_{i}$
$$
{\cal N}_{R_{0}/2}(f_{\lambda})\leq\mu_{f}(K)\ .
$$
These inequalities imply
$$
{\cal N}_{r}(f_{\lambda_{0}})=\mu_{f}(K)\ .
$$

The proof of the theorem is complete.\ \ \ \ \ \ $\Box$\\

\end{document}